\documentclass[10pt,twoside]{article}
\usepackage{Latex-document}

\def\bsh{\backslash}
\def\ify{\infty}
\def\ov{\overline}
\def\sbs{\subset}
\def\sps{\supset}
\def\ts{\times}

\def\a{\alpha}
\def\d{\delta}
\def\g{\gamma}
\def\om{\omega}
\def\vp{\varphi}

\def\D{\Delta}
\def\G{\Gamma}
\def\Om{\Omega}

\def\mpo{\mapsto}
\def\ra{\rightarrow}

\def\un{{\rm 1\mkern-4mu I}}

\def\Nc{{\cal N}}
\def\Oc{{\cal O}}

\font\tenbb=msbm10
\font\sevenbb=msbm7
\font\fivebb=msbm5
\newfam\bbfam
\textfont\bbfam=\tenbb
\scriptfont\bbfam=\sevenbb
\scriptscriptfont\bbfam=\fivebb
\def\bb{\fam\bbfam}

\def\Cb{{\bb C}}
\def\Nb{{\bb N}}
\def\Pb{{\bb P}}
\def\Qb{{\bb Q}}
\def\Rb{{\bb R}}

\def\build#1_#2^#3{\mathrel{
\mathop{\kern 0pt#1}\limits_{#2}^{#3}}}

\def\volumeno{I}

\title{\bf \boldmath Integrating $\partial \ov \partial$\vskip 6mm}
\author{Michael McQuillan\thanks{Institut des Hautes Etudes
Scientifiques, 35 route de Chartres, 91440 Bures-sur-Yvette,
France and Department of Mathematics, University of Glasgow
G128QW, Scotland. E-mail: mquillan@ihes.fr}\vspace{-0.5cm}}
\date{\vspace{-8mm}}

\begin{document}

\maketitle

\thispagestyle{first} \setcounter{page}{547}

\begin{abstract}

\vskip 3mm

We consider the algebro-geometric consequences of integration by
parts.

\vskip 4.5mm


\end{abstract}

\vskip 12mm

\section{Jensen's formula}

\vskip-5mm \hspace{5mm}

Recall that for a suitably regular function $\vp$ on the unit disc $\D$
we can apply integration by parts/Stoke's formula twice to obtain for
$r < 1$,
$$
\int_0^r \frac{dt}{t} \int_{\D (t)} dd^c \vp = \int_{\partial \D (r)}
\vp - \vp (0) \eqno (1.1)
$$
where $d^c = \frac{1}{4\pi i} \, (\partial - \ov \partial)$ so actually
we're integrating $\frac{1}{2 \pi i} \, \partial \ov \partial$.
In the presence of singularities things continue to work. For example
suppose $f : \D \ra X$ is a holomorphic map of complex spaces and $\ov
D$ a metricised effective Cartier divisor on $X$, with $f(0) \notin D$,
and $\vp = - \log f^* \Vert \un_D \Vert$, where $\un_D \in \Oc_X (D)$
is the tautological section, then we obtain,
$$
\matrix{
\displaystyle \int_0^r \frac{dt}{t} \int_{\D (t)} f^* c_1 (\ov D) &=
&\displaystyle - \int_{\partial \D (r)} \log \Vert f^* \un_D \Vert +
\log \Vert f^* \un_D \Vert (0) \cr \cr \cr
&+ &\displaystyle \sum_{0 < \vert z \vert < r} {\rm ord}_z (f^* D) \log
\left\vert \frac{r}{z} \right\vert \, . \hfill \cr
}
\eqno (1.2)
$$

Obviously it's not difficult to write down similar formulae for
not necessarily effective Cartier divisors, meromorphic functions,
drop the condition that $f(0) \notin D$ provided $f(\D) \not\sbs
D$, extend to ramified covers $p : Y \ra \D$, etc., but in all
cases what is clear is,

\medskip

\noindent {\bf Facts 1.3.} {\it (a) If $X$ is compact then $\int_0^r
\frac{dt}{t} \int_{\D (t)} f^* c_1 (\ov D)$ is very close to being
positive if $f(\D) \not\sbs D$, e.g. in the particular hypothesis
preceding (1.2),
$$
\int_0^r \frac{dt}{t} \int_{\D (t)} f^* c_1 (\ov D) \geq \log \Vert f^*
\un_D \Vert (0) + O_{\ov D} \, (1)
$$

(b) There is no such principle for the usual area function $\int_{\D
(r)} f^* c_1 (\ov D)$ except in extremely special cases such as $D$
ample.}

\medskip

Equally the essence of the study of curves on higher dimensional
varieties lies in understanding their intersection with divisors,
and, of course, the principle that a curve not lying in a divisor
intersects it positively is paramount to the discussion.
Consequently the right notion of intersection number for
non-compact curves is the so-called characteristic function
defined by either side of the identity (1.2). On the other hand
intersection number and integration are interchangeable in
algebraic geometry, and whence we will write,

\medskip

\noindent {\bf Notation 1.4.} {\it Let $\om$ be a $(1,1)$ form on $\D$ then
for $r < 1$,
$$
{\int\!\!\!\!\!\!\!\nabla}_{\!\!\D (r)} \om := \int_0^r \frac{dt}{t}
\int_{\D (t)} \om \, .
$$ }

The more traditional characteristic function notation is reserved for
the current associated to a map, i.e.

\medskip

\noindent {\bf Definition 1.5.} {\it Let $f : \D \ra X$ be a map of
complex spaces then for $r < 1$, we define,}
$$
T_f (r) : A^{1,1} (X) \ra \Cb : \om \mpo
{\int\!\!\!\!\!\!\!\nabla}_{\!\!\D (r)} f^* \om \, .
$$

Evidently in many cases one works with forms which are not quite
smooth, so there are variations on the definition. In any case in order
to motivate our intersection formalism let us pause to consider,

\section{Convergence}

\vskip-5mm \hspace{5mm}

The basic theorem in the study of subvarieties of a projective variety
is Grothen\-dieck's existence and properness of the Hilbert scheme, or
if one prefers a sequence of subvarieties of bounded degree has a
convergent subsequence. Of course families of smooth curves do not in
general limit on smooth curves but rather semi-stable ones, and as such
we must necessarily understand convergence of discs in the sense of
Gromov [G], i.e.

\medskip

\noindent {\bf Definition 2.1.} {\it A disc with bubbles $\D^b$ is a
connected $1$-dimensional complex space with singularities at worst nodes
exactly one of whose components is a disc $\D$ and such that every
connected component of $\D^b \bsh \{ \D \bsh {\rm sing} (\D^b)\}$ is a
tree of smooth rational curves.}

\medskip

For $z \in {\rm sing} (\D^b) \cap \D$ and $R_z$ the corresponding tree
of rational curves, and provided $0\notin {\rm sing} (\D^b)$ we can
extend our integral (1.4) to this more general situation by way of,
$$
{\int\!\!\!\!\!\!\!\nabla}_{\!\!\D^b (r)} \om :=
{\int\!\!\!\!\!\!\!\nabla}_{\!\!\D (r)} \om + \sum_{z \in \D (r) \cap
{\rm sing} (\D^b)} \log \frac{r}{\vert z \vert} \int_{R_z} \om
$$
while if $f : \D^b \ra X$ is a map then we have a graph,
$$
\G_f := ({\rm id} \ts f) (\D) \bigcup_{z \in \D \cap {\rm sing} (\D^b)}
z \ts f(R_z) \sbs \D \ts X \, .
$$

An appropriate formulation of Gromov's compactness theorem is then,

\medskip

\noindent {\bf Fact 2.2.} {\it Let $\ov{\rm Hom} (\D , X)$ be the space of
maps from discs with bubbles into a projective variety $X$ topologised
by way of the Hausdorff metric on the graphs then for $C : (0,1) \ra
\Rb_+$ any function and $K \sbs {\rm Aut} (\D)$ compact the set,
$$
\left\{ f \in \ov{\rm Hom} (\D , X) : \exists \, \a \in K \, , \
{\int\!\!\!\!\!\!\!\nabla}_{\!\!\D (r)} \a^* f^* c_1 (\ov H) \leq C(r)
\right\}
$$
where $\ov H$ is a metricised ample divisor, is compact.}

\medskip

Under more general hypothesis on $X$, 2.2 continues to hold, but in the
special projective case one has an essentially trivial proof thanks to
the ubiquitous Jensen formula, cf. [M]~V.3.1. Equally although the
appearance of automorphisms looks like an unwarranted complication they
are necessitated by,

\medskip

\noindent {\bf Remarks 2.3.} {(a) The possibility of bubbling at
the origin.

\smallskip

(b) The defect of positivity for the intersection number as per
(1.3)(a).}

\medskip

Observe moreover that the introduction of $\ov{\rm Hom} (\D , X)$ and
its precise relation to ${\rm Hom} (\D , X)$ are both necessary, and
easy respectively, i.e.

\medskip

\noindent {\bf Fact 2.4.} ([M] V.3.5) {\it  Let $T \sps {\rm Hom} (\D , X)$
be such that the bounded subsets in the sense of $2.2$ are relatively
compact then $T \sps \ov{\rm Hom} (\D , X)$. Moreover assuming that
$X$ is not absurdly singular then ${\rm Hom} (\D , X) = \ov{\rm Hom}
(\D , X)$ iff $X$ contains no rational curves.}

\medskip

One can equally generalise this to a log, or quasi-projective situation
by introducing a divisor $D$, whose components $D_i$ should be
$\Qb$-Cartier at which point the appropriate variation thanks to a
lemma of Mark Green, [Gr], is,

\medskip

\noindent {\bf Fact 2.5.} {\it  $\ov{\rm Hom} (\D , X \bsh D) \sbs {\rm Hom}
(\D , X)$ iff $X \bsh D$ and $D_i \bsh \build{\bigcup}_{j \ne i}^{}
D_j$ do not contain any affine lines.}

\medskip

In particular ${\rm Hom} (\D , X \bsh D)$ is relatively compact in
${\rm Hom} (\D , X)$ if and only if $\ov{\rm Hom} (\D , X \bsh D)$ is
compact and the boundary is {\it mildly hyperbolic} in the sense that $D_i
\bsh \build{\bigcup}_{j \ne i}^{} D_j$ does not contain affine lines.
The latter question is purely algebraic and closely related to the log
minimal model programme. In the case of foliations by curves an even
more delicate result holds since as Brunella has observed, [B], the
equivalence of $\ov{\rm Hom}$ with Hom for invariant maps into the
orbifold smooth part of a foliated variety is itself equivalent to the
said foliated variety being a minimal model.

\section{The Bloch principle}

\vskip-5mm \hspace{5mm}

Bloch's famous dictum, ``Nihil est in infinito quod non fuerit prius in
finito'', might thus be translated as,

\medskip

\noindent {\bf Question 3.1.} {Suppose for a projective variety
$X$, or more generally a log variety $(X,D)$ there is a Zariski
subset $Z$ of $X \bsh D$ through which every non-trivial map $f :
\Cb \ra X \bsh D$ must factor then do we have {\it hyperbolicity
modulo} $Z$, i.e. is it the case that a sequence $f_n$ in ${\rm
Hom} (\D , X)$ not affording a convergent subsequence in $\ov{\rm
Hom}$ must be arbitrarily close (in the compact open sense) to $Z
\cup D$.}

\medskip

In the particular case that 2.5 is satisfied we can replace $Z \cup D$
by $\ov Z$ and ask for {\it complete hyperbolicity modulo} $Z$, but
outside of surfaces (2.5) seems difficult to guarantee. Regardless in
his thesis Brody, [Br], provided an affirmative answer for both $Z$ and
$D$ empty by way of his reparameterisation lemma which was subsequently
extended by Green to the case of $Z$ empty and every $D_i \bsh
\build\bigcup_{j \ne i}^{} D_j$ not containing holomorphic lines.

Bearing in mind the singular variant of Green's lemma implicit in
2.5, which for example makes it applicable to stable families of
curves, it would appear that the unique known case not covered by
the methods of Brody and Green was a theorem of Bloch himself,
[Bl], i.e. $\Pb^2 \bsh \{$4 planes in general position$\}$, and
its subsequent extension by Cartan to $\Pb^n$, [C]. However, even
here, a moment's inspection shows that 2.5 holds, so one knows a
priori that there can be no bubbling, and whence complete
hyperbolicity in the sense of 3.1 trivially implies so-called
normal convergence modulo the diagonal hyperplanes, and the
correct structure is obscured.

\medskip

Now an extension of the reparameterisation lemma to cover 3.1 would be
by far the most preferable way forward, since the non-existence of
holomorphic lines is an essentially useless qualitative statement
without the quantitative information provided by the convergence of
discs. Nevertheless we can vaguely approximate a reparameterisation
lemma thanks as ever to Jensen's formula. Specifically consider as
given,

\medskip

\noindent {\bf Data 3.2.}{\it \begin{itemize}

\item[(a)] A $\Qb$-Cartier divisor $\partial$ on a
log-variety $(X,D)$.

\item[(b)] A sequence $f_n \in {\rm Hom} (\D , X \bsh D)$ which neither
affords a convergent subsequence nor is arbitrarily close to
$\partial \cup D$.
\end{itemize} }

In light of (b) we can choose convergent automorphisms $\a_n \in
{\rm Aut} (\D)$, such that $\a_n^* f_n (0)$ is bounded away from
$\partial \cup D$, and given, modulo subsequencing, the
convergence of the $\a_n$ we may as well suppose this. Moreover
for each $0 < r < 1$ we can normalise the current $T_{f_n} (r)$ of
1.5 by its degree with respect to an ample divisor $H$, which
we'll denote by $T_{f_n}^H (r)$ and take a weak limit for a
suitable subset $\Nc$ of $\Nb$ to obtain a current $T_{\Nc}^H
(r)$. In addition 3.2(b) also tells us that for some fixed $0 < s
< 1$, the degrees of the $f_n$ at $s$ go to infinity, and whence
by (1.1) and (1.2)

\medskip

\noindent {\bf Pre-Fact 3.3.} {\it  For $r \geq s$, $T_{\Nc}^H (r)$ is a
positive harmonic current such that, $T_{\Nc}^H (r) \cdot F \geq 0$ for
all effective divisors $F$ supported in $\partial \cup D$.}

\medskip

What is somewhat less trivial, but once more the key is Jensen's
formula, is,

\medskip

\noindent {\bf Fact 3.3(bis).} ([M] V.2.4) {\it Subsequencing in $\Nc$ as
necessary, then for $r \geq s$ outside of a set of finite hyperbolic
measure (i.e. $(1-r^2)^{-1} dr$) $T_{\Nc}^H (r)$ is closed.}

\medskip

Obviously there are various choices involved but whenever we're
dealing in the context of countably many projective varieties they
can all be rendered functorial, up to a constant, with respect to
push forward. The constant itself only causes a problem should it
be zero which is usually what one wants to prove anyway, and as
such the notation $T(r)$ is relatively unambiguous, and represents
in a vague sense a parabolic limit of the sequence $f_n$.

\section{Applications}

\vskip-5mm \hspace{5mm}

Applications of course require some knowledge of intersection numbers,
and quite generally even for a compact curve $f : C \ra X$ there is
very little that one can say in general beyond,

\medskip

\noindent {\bf Observation 4.1.} {\it Let $f' :C \ra P(T_X)$ be the
derivative $(\Pb (\Om_X)$ in the notation of {\rm EGA)} with $L$ the
tautological bundle then,
$$
L_{\cdot_{f'}} C = (2g-2) - {\rm Ram}_f \, .
$$ }

This is of course the Riemann-Hurwitz formula if $\dim  X =1$, and
there's an equally trivially log-variant where on the right hand
side we have to throw in the number of points in the intersection
with the boundary $D$ counted without multiplicity the special
case of $\Pb^1 \bsh \{ 0 , 1 , \ify \}$ being Mason's ``$a,b,c$''
theorem for polynomials. The correct generality for best possible
applications is to work with log-smooth Deligne-Mumford stacks (or
alternatively just orbifolds since the inertia tends to be
irrelevant), however for simplicity let's stick with log-smooth
varieties and metricise $T_X (-\log D)$ by way of a complete
metric $\Vert \ \Vert_{\log}$ on $(X,D)$, which in turn leads to a
mildly singular metricisation $\ov L$ of the tautological bundle.
Supposing for simplicity that $f(0) \notin D$ with $f$ unramified
at the origin then Jensen's formula yields,

\medskip

\noindent {\bf Observation 4.2(bis).} {\it  Notations as above,}
$$
\matrix{
\displaystyle {\int\!\!\!\!\!\!\!\nabla}_{\!\!\D (r)} f^* c_1 (\ov L)
&= &\displaystyle -\log \left\Vert f_* \left( \frac{\partial}{\partial
z} \right) \right\Vert_{\log} (0) + \int_{\partial \D(r)} \log
\left\Vert f_* \left( \frac{\partial}{\partial z} \right)
\right\Vert_{\log} \hfill \cr \cr \cr
&+ &\displaystyle \sum_{0 < \vert z \vert < r} \min \{ 1 , {\rm ord}_z
(f^* D)\} \log \frac{r}{\vert z \vert} \hfill \cr \cr \cr
&- &\displaystyle \sum_{{0 < \vert z \vert < r \atop f(z) \notin D}}
{\rm ord}_z (R_f) \log \frac{r}{\vert z \vert} \, . \hfill \cr
}
$$

Combining the concavity of the logarithm and once more Jensen's formula,
but this time for $dd^c \log \log^2 \Vert \un_D \Vert$ for any norm on
the boundary divisor $D$, immediately yields in the notations of 3.3,

\medskip

\noindent {\bf Fact 4.3.} {\it  Let $T'(r)$ be the current associated to the
logarithmic derivative of a sequence $f_n \in {\rm Hom} (\D , X \bsh
D)$ with $f_n (0)$ not arbitrarily close to $D$ and which does not
afford a convergent subsequence then outwith a set of finite hyperbolic
measure,}
$$
L_{\cdot} \, T'(r) \leq 0 \, .
$$

The so-called tautological inequality 4.3 is well adapted for
applications to convergence of discs (note incidentally that it's
implicit to the formulation that a smooth metric on the bundle
$T_X (-\log D)$ is being employed). Nevertheless for more delicate
questions such as quantifying degenerate/non-convergent behaviour.
etc. there is a wealth of information in (4.2) that is lost in the
coarser corollary. Indeed even using the concavity of the
logarithm distorts a very delicate term measuring the
`ramification at $\infty$', i.e. the distorsion of the boundary
from it's length in the Poincar\'e metric, which is closely
related to the difficulty of extracting an isoperimetric
inequality from a knowledge of hyperbolicity in the sense of 3.1.
While from the still deeper curvature point of view, 4.2 is simply
a doubly integrated tautological Schwarz lemma, since by
definition metricising $T_X(-\log D)$ by way of a metric $\omega$
of curvature $\leq -K$ is equivalent to a lower bound of the left
hand side of the form,

$$K\,\,
{\int\!\!\!\!\!\!\!\nabla}_{\!\!\D (r)} \om $$

\noindent for all infinitely small, and whence all in the large, possible discs.
While on the subject of curvature and
isoperimetric inequalities, a variant specific to
dimension 1 replaces the current $\d_D$ implicitly hidden in (4.2) by
the current associated to the boundary of a simply connected region,
i.e.

\medskip

\noindent {\bf Variant 4.4.} {\it Suppose $\dim X = 1$, let $U_i \sbs X$ be
simply connected, $h_i : \D \build\longrightarrow_{}^{\sim} U_i$
isomorphisms and put
$$
\d_{\G_i} : A^0 (X) \ra \Cb : \vp \mpo \int_{\partial \D} h_i^* \vp \,
.
$$
Then specific to dimension 1, $\d_{\G_i}$ is closed and may be
written, $c_1 (H) + dd^c \g_i$ for $H$ an ample divisor of degree 1.
Now apply Jensen's formula to recover an integrated form of Ahlfor's
isoperimetric inequality, and the Five Island's Theorem.}

\medskip

Returning to varieties and divisors it's still possible to employ (4.2)
to get integrated isoperimetric inequalities for more general
situations that preserve some 1-dimensional flavour, i.e. discs which
are invariant by foliations by curves, with {\it canonical foliation
singularities} (with the obvious definition of that notion which is
functorial with respect to the ideas) and which do not pass through the
singularities. The latter hypothesis which is reasonable for the study
of the leafwise variation of the Poincar\'e metric is however somewhat
restrictive for other applications, and is probably unnecessary as
suggested by the essentially optimal inequality of [M]~V.4.4 for
foliations on surfaces which employs (4.2) to a very large number of
monoidal transformations in the foliation singularities. Regardless
here is a genuinely 2-dimensional theorem,

\medskip

\noindent {\bf Theorem 4.5.} ([M] V.5) {\it Let $(X,D)$ be a smooth
logarithmic surface with $\Om_X (\log D)$ big (e.g. log-general type,
and $s_2 (\Om_X (\log D)) > 0$) then there is a proper Zariski subset
$Z$ of $X \bsh D$ such that $X \bsh D$ is complete hyperbolic (in the
sense of $3.1$ et sequel) modulo $Z$.}

Indeed one can even optimally quantify (cf. op. cit.) the degeneration
of the Kobayashi metric (which is evidently continuous and non-zero off
$Z$) around $Z$. Amusingly the theorem only covers $\Pb^2 \bsh \{$5
planes in general position$\}$, although it's a good exercise in the
techniques (cf. op. cit. V.4) to prove Bloch's theorem too, at which
point a rather small sequence of blow ups replaces all of the original
estimation. In any case (4.5) should only be seen as a stepping stone
which in order of ascending difficulty leaves open the following
questions, viz,

\medskip

\noindent {\bf Concluding Remarks 4.6.} { For concreteness take a
smooth algebraic surface $X$ of general type with $c_1^2 > c_2$
(otherwise the following should be understood in terms of higher
jets, but not for anything more general than a surface) then,
\begin{itemize}
\item[(a)] Do we have an isoperimetric inequality with appropriate
degeneration along the subset $Z$ of (4.5).
\item[(b)] Is the Kobayashi metric negatively curved.
\item[(c)] For each $x \notin Z$ and $t$ a tangent direction at $x$, is
there a unique up to the usual action of ${\rm SL}_2 (\Rb)$ pointed
disc with maximal tangent
in the direction $t$, and if so does it continue to be so along
its image, i.e. is there a continuous (off $Z$) connection whose
geodesics are the discs defining the Kobayashi metric.
\end{itemize}
}

\section{Thanks}

\vskip-5mm \hspace{5mm}

In closing it
is a pleasure to thank M.~Brunella for introducing me to the unit
disc and bubbling, together with M.~Gromov for continuing
this education, but above all C\'ecile without whom the reader could
never have got this far.

\label{lastpage}

\end{document}